\documentclass[a4paper, parskip=full+, bibliography=totocnumbered,
leqno,
headsepline, %
twoside, openright
]{scrartcl}
\pdfoutput=1
\usepackage{scrpage2}

\usepackage[OptionalSectionArgument=Head]{screxperimental}
\usepackage{lipsum}

\lehead{T. Wieber}
\rehead{\titelkurz}

\lohead{\rightmark}
\rohead{}

\automark[section]{chapter}

\usepackage[dvipsnames  
]{xcolor}

\usepackage[normalem]{ulem}

\usepackage[T1]{fontenc}
\usepackage[latin1]{inputenc}

\usepackage{relsize}

\usepackage{amsmath}
\usepackage{amssymb}
\usepackage{amsmath}

\allowdisplaybreaks[1]

\usepackage[amsmath,thmmarks,hyperref]{ntheorem}
\usepackage{varioref}

\usepackage{amsgen}
\usepackage{amscd}
\usepackage{nicefrac}
\usepackage{extarrows}

\usepackage{listings}

\usepackage{tikz}
\usetikzlibrary{matrix,arrows,intersections,shapes.arrows,fit}
\usepackage{epstopdf}

\usepackage[colorinlistoftodos]{todonotes}

\usepackage[refpage,intoc]{nomencl}

\usepackage{makeidx}

\makeindex

\makenomenclature

\usepackage[
	backref=page,
	pdfauthor={Thomas Wieber},
	pdfsubject={We develop two structure theorems for vector valued Siegel modular forms for Igusa's subgroup Gamma_2[2,4], the multiplier system induced by the theta constants and the representation Sym^2. In the proof, we identify some of these modular forms with rational tensors with easily handleable poles on P^3C. It follows that the observed modules of modular forms are generated by the Rankin-Cohen brackets of the four theta series of the second kind.},
	pdfcreator={AUCTeX, LaTeX with hyperref and KOMA-Script},
	pdfkeywords={vector valued Siegel modular form, Hilbert
          function, algebraic geometry, complex differential geometry,
          tensor field, generators, Igusa's subgroup of level 2 Gamma _2[2, 4], structure theorem, symmetric square of the standard representation Sym2, meromorphic tensor},
	pdfpagemode=UseOutlines,
	pdfdisplaydoctitle=true,
	pdflang=eng
]{hyperref}

   \renewcommand*{\backref}[1]{%
   }%
   \renewcommand*{\backrefalt}[4]{%
   	\mbox{(cited on %
   	\ifnum#1=1 %
		   page~%
	   \else
   		pages~%
   	\fi
   	#2)}%
   }

\makeatletter
\def\thenomenclature{%
  \@ifundefined{addsec}%
  {
    \@ifundefined{chapter}%
    {
      \section*{\nomname}
      \if@intoc\addcontentsline{toc}{section}{\nomname}\fi%
    }%
    {
      \chapter*{\nomname}
      \if@intoc\addcontentsline{toc}{chapter}{\nomname}\fi%
    }%
  }%
  {
    \@ifundefined{addchap}%
    {
      \if@intoc\addsec{\nomname}\else\addsec*{\nomname}\fi%
    }%
    {
      \if@intoc\addchap{\nomname}\else\addchap*{\nomname}\fi%
    }%
  }%
  \nompreamble
  \list{}{%
    \labelwidth\nom@tempdim
    \leftmargin\labelwidth
    \advance\leftmargin\labelsep
    \itemsep\nomitemsep
    \let\makelabel\nomlabel}}
\makeatother

\usepackage{cleveref}

\crefname{cor}{corollary}{corollaries}
\Crefname{cor}{Corollary}{Corollaries}
\crefname{lem}{lemma}{lemmas}
\Crefname{lem}{Lemma}{Lemmas}
\crefname{rem}{remark}{remarks}
\Crefname{rem}{Remark}{Remarks}
\crefname{defi}{definition}{definitions}
\Crefname{defi}{Definition}{Definitions}
\crefname{thm}{theorem}{theorems}
\Crefname{thm}{Theorem}{Theorems}
\crefname{cex}{counterexample}{counterexamples}
\Crefname{cex}{Counterexample}{Counterexamples}
\crefname{ex}{example}{examples}
\Crefname{ex}{Example}{Examples}
\crefname{conj}{conjecture}{conjectures}
\Crefname{conj}{Conjecture}{Conjectures}

\newenvironment{eq*}{\begin{equation*}}{\end{equation*}}

\newcommand{\Adj}{\ensuremath{\operatorname{Adj}}}

\renewcommand{\Im}{\operatorname{Im}}

\newcommand{\supp}{\operatorname{supp}}

\newcommand{\ord}[1]{\ensuremath{\operatorname{ord}\left( #1 \right)}}

\newcommand{\Aut}[1][V]{\ensuremath{
\operatorname{Aut}\left( { #1 } \right)}}

\newcommand{\Altk}[2][k]{\ensuremath{\Lambda^{#1}#2}}

\newcommand{\Fix}[1][\G]{\ensuremath{\operatorname{S}\left( #1 \right)}}

\newcommand{\ohne}{\setminus}

\newcommand{\iso}{\cong}

\newcommand{\union}{\bigcup}

\newcommand{\we}{\wedge}
\newcommand{\del}{\ensuremath{\partial}}

\newcommand{\abs}[1]{\ensuremath{\left| #1 \right|}}
\newcommand{\set}[1]{\ensuremath{\left\{ #1 \right\}}}
\newcommand{\eck}[1]{\ensuremath{\left[ #1 \right]}}

\newcommand{\floor}[1]{\ensuremath{\left\lfloor #1 \right\rfloor}}

\newcommand{\ten}{\otimes}

\newcommand{\Tenk}{\bigotimes_{l=1}^{k}} 

\newcommand{\tenO}{\otimes_{\Hol}}

\newcommand{\vt}{\ensuremath{\vartheta}}

\newcommand{\nui}{\ensuremath{{\nu^i}}}

\newcommand{\p}{\ensuremath{\pi}}

\newcommand{\rh}{\ensuremath{\rho}}

\newcommand{\om}{\ensuremath{\omega}}

\newcommand{\G}{\ensuremath{\Gamma}}
\newcommand{\Gt}{\ensuremath{\Gamma_2}}
\newcommand{\Gn}[1][n]{\ensuremath{{\Gamma_{#1}}}}
\newcommand{\Gnq}{\ensuremath{{\Gamma_n\left[q \right]}}}
\newcommand{\Gttf}{\ensuremath{{\Gamma_2\left[ 2,4 \right]}}}
\newcommand{\Gttfpm}{\ensuremath{\left(\nicefrac{\Gttf}{\pm I_4}\right)}}

\newcommand{\Om}{\ensuremath{\Omega}}

\newcommand{\entcov}{\ent\ covering}
\newcommand{\anentcov}{a simple covering}
\newcommand{\ent}{simple}
\newcommand{\fn}{function}
\newcommand{\hol}{holomorphic}
\newcommand{\holfn}{holomorphic function}
\newcommand{\holten}{holomorphic tensor}

\newcommand{\meroten}{meromorphic tensor}

\newcommand{\Obda}{Without loss of generality}
\newcommand{\nec}{necessary}

\newcommand{\nhd}{neighbourhood}
\newcommand{\resp}{respectively}

\newcommand{\am}[1][\mchar]{\ensuremath{a_{#1}}}

\newcommand{\magma}{\textsc{magma}\cite{magma}}

\newcommand{\zt}{\ensuremath{\mathfrak{z}}}

\newcommand{\dzti}[1][i]{\ensuremath{d\mathfrak{z}^{#1}}}

\newcommand{\dzt}{\ensuremath{d\mathfrak{z}}}

\newcommand{\fZtail}{
\ensuremath{f(Z) \cdot
\left(\dZ^0 \we \dZ^1 \we \dZ^2 \right)^{\ten \kforl}}
}

\newcommand{\dzbig}{\ensuremath{dZ}}
\newcommand{\dZ}{\dzbig}
\newcommand{\dZi}[1][i]{\ensuremath{\dzbig^{#1}}}
\newcommand{\dzlex}{\ensuremath{\bigwedge\dZij}}

\newcommand{\dzij}[1][ij]{\ensuremath{\dZ ^{#1}}}
\newcommand{\dZij}[1][ij]{\dzij[#1]}

\newcommand{\delzti}[2][i]{\ensuremath{\frac{\partial #2}{\partial \zt^#1}}}
\newcommand{\delztj}[1]{\delzti[j]{#1}}

\newcommand{\homproj}[1]{\ensuremath{\pi_{*}\left(#1\right)}}

\newcommand{\diff}{\ensuremath{d}}

\newcommand{\rhc}{\ensuremath{\rho_e^*}}
\newcommand{\rhe}{\ensuremath{\rho_e}}
\newcommand{\rhdmz}{\ensuremath{\rhc\! \ten\! \rhc}}

\newcommand{\rhdmzsym}[1][n]{
\ensuremath{ \rhc\! \odot\! \rhc}}

\newcommand{\rhst}{\ensuremath{\rhe\! \ten\! \rhe}}

\newcommand{\rhstsym}[1][n]{\ensuremath{\operatorname{Sym}^2}}
\newcommand{\rhsym}{\ensuremath{\operatorname{Sym}^2}} 

\newcommand{\JtenJ}{\ensuremath{J^{-1}_1\! \ten\! J^{-1}_1}}
\newcommand{\JtenJsymt}{\ensuremath{ \bound{\JtenJ}{\sym[2]}}}

\newcommand{\J}{\ensuremath{{\mathcal{J}}}}
\newcommand{\Jr}[1][r]{\ensuremath{\mathfrak{J}^{#1}}}

\newcommand{\sqrtr}[2][r]{\ensuremath{\sqrt{#2}\;\strut^{#1} }}

\def\standardmultiplierpower{1}
\newcommand{\multipliersystem}[2]{%
  \begingroup
    \def\tmp{#2}%
    \ifx\tmp\standardmultiplierpower
        \ensuremath{#1}%
    \else
        \ensuremath{#1^{#2}}%
    \fi
  \endgroup}

\newcommand{\vf}[1][1]{\multipliersystem{v_{f}}{#1}}

\newcommand{\thesermZ}{\ensuremath{\vartheta\!\left[\!\begin{array}{c}\mchar^1
\\ \mchar^2 \end{array}\!\right]\!(\!Z\!)}}
\newcommand{\theserqmZ}{\ensuremath{\vartheta^2\!\left[\!\begin{array}{c}
\mchar^1 \\ \mchar^2 \end{array}\!\right]\!(\!Z\!)}}
\newcommand{\chifive}{\ensuremath{\chi_5}}

\newcommand{\ramdiv}{\ensuremath{{\cal R}}}

\newcommand{\mfix}{\ensuremath{\mathcal{M}}}
\newcommand{\mchar}{\ensuremath{\mathfrak{m}}}
\newcommand{\mcharanz}{\ensuremath{M}}

\newcommand{\Pe}[1][i]{\ensuremath{\mathsf{P}_{#1}}}

\newcommand{\Qm}[1][\mchar]{\ensuremath{Q_{#1}}}

\newcommand{\R}{\ensuremath{\mathbb{R}}}

\newcommand{\C}{\ensuremath{\mathbb{C}}}

\newcommand{\Cn}[1][n]{\ensuremath{\mathbb{C}^{#1}}}

\newcommand{\N}{\ensuremath{\mathbb{N}}}
\newcommand{\Z}{\ensuremath{\mathbb{Z}}}
\newcommand{\Zq}[1][q]{
\ensuremath{\nicefrac{\mathbb{Z}}{#1\mathbb{Z}}}}
\newcommand{\Zn}[1][n]{\ensuremath{\mathbb{Z}^{#1}}}

\newcommand{\F}{\ensuremath{\mathbb{F}}}

\newcommand{\E}{\ensuremath{\mathbb{E}}} 
\newcommand{\En}[1][n]{\ensuremath{\mathbb{E}^{#1}}}

\newcommand{\PnC}[1][n]{\ensuremath{\mathbb{P}^{#1}\mathbb{C}}}
\newcommand{\A}{\ensuremath{\mathbb{A}}}
\newcommand{\Anull}{\ensuremath{\mathbb{A}_{0}}}

\newcommand{\Jn}[1][n]{\ensuremath{J_{#1}}}

\newcommand{\GL}{\operatorname{GL}}
\newcommand{\GLnC}[1][n]{\ensuremath{\GL(#1,\C)}}

\newcommand{\SL}{\operatorname{SL}}

\newcommand{\Sp}{\operatorname{Sp}}

\newcommand{\SpnZ}[1][n]{\ensuremath{\Sp(#1,\Z)}}
\newcommand{\SpnZq}[1][n]{\ensuremath{\Sp(#1,\Zq)}}
\newcommand{\SpnR}[1][n]{\ensuremath{\Sp(#1,\R)}}
\newcommand{\Spnr}[1][n]{\ensuremath{\Sp(#1,R)}}

\newcommand{\sym}[1][n]{\ensuremath{
\operatorname{Sym}^2\!\left(\Cn[#1]\right)}}

\newcommand{\Hn}[1][n]{\ensuremath{\mathbb{H}_{#1}}}
\newcommand{\Ht}{\ensuremath{\mathbb{H}_{2}}}
\newcommand{\Dn}[1][n]{\ensuremath{\Delta_{#1}}}

\newcommand{\Hol}{\ensuremath{{\cal O}}}
\newcommand{\OU}[1][U]{\ensuremath{{\cal O}(#1)}} %

\newcommand{\Omp}[1][p]{\ensuremath{ \Altk[#1]{\Omega}}}
\newcommand{\Omn}{\Omp[n]}

\newcommand{\Oot}[1][n]{\ensuremath{
\left(\Om \tenO \left(\Omp[#1]\right)^{\ten k}\right)}}

\newcommand{\Opt}[1][n]{\ensuremath{\Omp \tenO (\Omp[#1] )^{\ten k}}}

\newcommand{\Ott}[1][n]{\ensuremath{
\left(\Omp[2]\tenO \left(\Omp[#1]\right)^{\ten k}\right)}}

\newcommand{\Oiso}{\ensuremath{
\left(\Omp[(\frac{n(n+1)}{2}-1)]\tenO \left(\Omp[\frac{n(n+1)}{2}]\right)^{\ten k}\right)}}

\newcommand{\Ot}[1][k]{\ensuremath{
\Omega
^{\ten #1}}}

\newcommand{\kforl}{\ensuremath{r}}
\newcommand{\kforls}{\ensuremath{s}}

\def\standardmultiplierweight{1}
\newcommand{\multiplierweight}[1]{%
  \begingroup
    \def\tmp{#1}%
    \ifx\tmp\standardmultiplierweight
        \ensuremath{\frac{1}{2}}%
    \else
        \ensuremath{\frac{\tmp}{2}}%
    \fi
  \endgroup}

\newcommand{\scamodfunth}[1][1]{
\ensuremath{\left[ \Gttf ,\frac{#1}{2}
,\vf[#1] \right]}}

\newcommand{\scamodfunthvf}[1][1]{
\ensuremath{\left[ \Gttf ,\frac{#1}{2},\vf[#1] \cdot \vf[2] \right]}}
\newcommand{\aoki}[3][]{\ensuremath{\left\{ #2,  #3  \right\}_{#1}}}
\newcommand{\rcb}[4][]{\ensuremath{\left\{ #2,  #3,  #4  \right\}_{#1}}}

\newcommand{\Attf}[1][]{\ensuremath{A^{int}_{#1}(\Gttf)}}
\newcommand{\Attfvf}[1][]{\ensuremath{A_{#1}(\Gttf,\vf)}} 

\newcommand{\Attfrhvf}[1][]{\ensuremath{{\cal M}^+_{#1} (\Gttf)}}
\newcommand{\Attfrhvfvf}[1][]{\ensuremath{{\cal M}^-_{#1} (\Gttf)}}

\newcommand{\Mp}[1][]{\ensuremath{{\cal M}^+_{#1}}}
\newcommand{\Mm}[1][]{\ensuremath{{\cal M}^-_{#1}}}
\newcommand{\Mmrcb}[1][]{\ensuremath{M^-_{#1}}}

\newcommand{\vecmodfun}[1][]{\ensuremath{
\left[ \Gamma,\left( \frac{r}{2},\rh \right),\multipliersystem{v}{#1} \right]}}
\newcommand{\vecmodfunst}[1][(n+1)k]{\ensuremath{
\left[ \Gamma,\left( #1,\rhstsym \right)\right]}}
\newcommand{\vecmodfunstsymt}[1][3k]{\ensuremath{
\left[ \Gamma, \left(#1,\rhstsym[2] \right)\right]}}
\newcommand{\vecmodfunstpm}[1][3\kforl]{\ensuremath{
\left[ \Gttf, \left( #1,\rhstsym[2] \right)\right]}}
\newcommand{\vecmodfunstvf}[1][r]{\ensuremath{
\left[ \Gttf,\left( \frac{#1}{2},\rhstsym\right),\vf[#1]\right]}}
\newcommand{\vecmodfunstvfvf}[1][r]{\ensuremath{
\left[ \Gttf,\left( \frac{#1}{2},\rhstsym\right),\vf[#1]\cdot \vf[2]\right]}}
\newcommand{\vecmodfundmzsymt}[1][3k]{\ensuremath{
\left[ \Gamma, \left(#1,\rhdmzsym[2] \right)\right]}}
\newcommand{\polyfa}[1][]{\ensuremath{\C_{#1}[f_0, \dots , f_3]}}
\newcommand{\polyfab}[1][]{\ensuremath{\left(\C_{#1}[f_0, \dots , f_3]\right)}}

\newcommand{\zero}[1]{\ensuremath{Z\left(#1\right)}}

\newcommand{\zeropi}[1]{\ensuremath{\pi\left(Z\left(#1\right)\right)}}

\newcommand{\pnk}[1][k]{\ensuremath{ p^{#1}_{n}}}

\newcommand{\rami}[2]{\ensuremath{#1\left({\operatorname{Ram}}\left(#1\right)\right)}}

\newcommand{\longhence}{\Longrightarrow}
\newcommand{\hooklongrightarrow}{\lhook\joinrel\longrightarrow}

\newcommand{\map}{\rightarrow}
\newcommand{\longmap}{\longrightarrow}
\newcommand{\longto}{\longrightarrow}

\newcommand{\longinjmap}{\hooklongrightarrow}

\newcommand{\inpro}[1]{\ensuremath{\left\langle #1 \right\rangle}}
\newcommand{\MZ}[1][Z]{\ensuremath{M\!\left\langle #1 \right\rangle}}

\newcommand{\DMZ}[1][Z]{\ensuremath{DM\!\left\langle #1\right\rangle}}

\newcommand{\lecos}[2]{\ensuremath{%
  \hbox{\kern 0.1em%
  \lower 0.25ex\hbox {\scriptsize$#2$}%
  \raise 0.25ex\hbox {\scriptsize$\backslash$}%
  \kern -0.1em%
  \raise 0.5ex\hbox {\scriptsize$#1$}}%
  \kern  0.2em}}

\newcommand{\inv}[1]{\ensuremath{\frac{1}{#1}}}

\newcommand{\bound}[2]{
\ensuremath{\left.\left( #1 \right) \right|_{#2}}}

\newcommand{\diag}[1]{\ensuremath{\left( #1 \right) _{0}}}

\newcommand{\kasten}[1]{
\fbox {
    \parbox{0.967\linewidth}{#1}
}
}

\newcommand{\paperempty}[5]{}

\newcommand{\abb}[5]
{%
\[
\begin{array}{lccl}
 #1 : & #2 &\longmap & #3, \\
 & #4 &\longmapsto & #5
\end{array}
\]
}

\newcommand{\pagsref}[2]{on pages \pageref{#1} and \pageref{#2}, \resp
}

\newcommand{\pref}[1]{\cref{#1} \vpageref{#1}}
\newcommand{\Pref}[1]{\Cref{#1} \vpageref{#1}}
\newcommand{\psref}[2]{\cref{#1,#2} \pagsref{#1}{#2}}
\newcommand{\prefs}[2]{\psref{#1}{#2}}

\newcommand{\skizze}[1]{}

\newcommand{\todocom}[2][]{\todo[noline,color=yellow!60,{#1}]{\thesubsection.
#2}}

\newcommand{\nc}[3][]{\nomenclature[#1]{#2}{#3}}

\newcommand{\beginroman}
{
\begin{enumerate}
\renewcommand{\theenumi}{\roman{enumi}}
\item
}

\newcommand{\beginromanchange}
{
\begin{enumerate}
\renewcommand{\theenumi}{\roman{enumi}}
}

\newcommand{\beginarabic}
{
\begin{enumerate}
\item
}

\newcommand{\titel}{Vector Valued 
Siegel Modular Forms \\
for \Gttf \ and \rhstsym }
\newcommand{\titelkurz}{ Vector Valued 
Siegel Modular Forms
for \Gttf \ and \rhstsym }
\title{\titel}

\date{}

\author{Thomas Wieber}

\begin{document}
\makeatletter
\g@addto@macro\normalsize{
\setlength\abovedisplayskip{0pt}
\setlength\belowdisplayskip{0pt}
}
\makeatother


\maketitle
\begin{abstract}\noindent
We develop two structure theorems for vector
valued Siegel modular forms for Igusa's subgroup \Gttf,
the multiplier system induced by the theta constants and the
 representation \rhsym . In the proof, we identify some of  these modular forms with rational tensors with easily handleable poles on~\PnC[3]. It follows that the observed modules of modular forms are generated by the
Rankin-Cohen brackets of the four theta series of the second kind.
\end{abstract}

\section{Introduction}
\label{chap:intro}

A vector valued modular form with respect to a subgroup $\G \subset
\SpnZ$ of finite index is a \hol\ function  $f:\Hn \map V$ that transforms
under \G \ as follows :%
$$f(\MZ)=j(M,Z)f(Z).$$
Here $j$ denotes a factor of
automorphy on a finite dimensional vector space $V$, i.e. a map
$j:\G\times\Hn \to \GL(V)$ that is \hol \ in
the second variable and satisfies the cocycle relation $j(MN,Z)=j(M,N\inpro{Z})j(N,Z)$.
\newline
The characters of \G\ are easy examples of factors of automorphy. Furthermore, a given polynomial representation $\rh : \GLnC \map
\GL(V)$ induces the factor of automorphy $\rh((CZ+D))$. For integral $r$, the functions
$\sqrtr{det(CZ+D)}$ satisfy the cocyle relation up to
$\pm 1$. In order to compensate this error, they are multiplied with
multiplier systems~$v(M)$ of weight $r/2$.

In \cite{tsushima}, Tsushima calculated the dimension of vector
spaces of vector valued cusp forms by means of the
Riemann-Roch-Hirzebruch-Theorem.
In \cite{satoh}, Satoh combined this result with the decomposition of the
vector space of vector valued modular forms into the subspace of cusp
forms and the subsapce of Eisenstein series, cf. \cite{arakawa}.
He obtained a structure theorem for certain vector spaces of
vector valued modular forms with respect to the full modular group \SpnZ[2].
We should also mention \cite{ibukiyama} and  \cite{dorp} which used the same strategy.
Similarly, \cite{aoki} and \cite{geer} estimated or calculated the Hilbert series and found matching generators.

In this article, we use a geometric method to get similar results for Igusa's group 
$$\G_2\eck{2,4}:=\set{M \in \G_2\eck{2} : \diag{AB^t}\equiv\diag{CD^t}\equiv0
\mod 4}$$
 instead of the full modular group. 
We benefit from the fact that the Satake compactification of
\lecos{\ \Ht}{\Gttf\ } is simply the 3 dimensional projective space
\PnC[3]. This is a consequence of a couple of basic results by Igusa \cite[Going-down process, p.397]{igusa} that was proven by Runge \cite{runge}.

We shall investigate modular forms being linked to the symmetric
square of the standard representation \rhe, i. e.
\[ \rhstsym=\operatorname{Sym}^2(\rhe):\GLnC[2]
\longto \Aut[
{\sym[2]}
].\]

%
Here \sym[2]\ denotes the symmetric square of \Cn[2], which can be identified with the
space of symmetric $2 \times 2$ matrices.
Then, the representation is given by
\[ {A} \longmapsto \set{X \mapsto AXA^{t}}.\]
%

In particular, we shall study the spaces of modular forms with respect to  representations of the type
\[
\det{} ^k \ten \left(\rhstsym \right), \quad\quad 
k \in \Z.
\]
For these vector spaces we shall give generators and compute the
dimensions. This can be found either on the following pages or in \pref{dimofvecmod}.

Our approach relies on the fact that in  the cases where $k=3r$ and $k=3r+1$, these forms
 can be identified with \Gttf - invariant \holten s of the form%
%
\begin{align}
\label{eventype}
& \left( f_{0}(Z) \cdot \dZi[0]+f_{1}(Z)  \cdot  \dZi[1]
+f_{2}(Z)  \cdot  \dZi[2] \right)\ten
\left(\dZ^0 \!\wedge \! \dZ^1 \!\wedge \! \dZ^2 \right)^{\ten \kforl} 
\intertext{and}
\label{oddtype}
& \left( g_{0}(Z)  \cdot  \dZi[1]\!\wedge \!\dZi[2]+g_{1}(Z)  \cdot
  \dZi[0]\!\wedge \!\dZi[2]
+g_{2}(Z)  \cdot  \dZi[0]\!\wedge \!\dZi[1] \right) 
\ten
\left(\dZ^0 \!\wedge \! \dZ^1 \!\wedge \! \dZ^2 \right)^{\ten \kforl} ,
\end{align}
\resp . Here the points in \Ht \ are of the form 
\[ Z=\begin{pmatrix}
Z_0 & Z_1 \\
Z_1 & Z_2
\end{pmatrix}
. \]

The crucial fact is that the map
\[ \Ht \longto 
\lecos{\Ht}{\Gttf}\longinjmap \PnC[3] \]
branches over 10 explicitly given quadrics in \PnC[3]. This implies
that \Gttf -invariant tensors on \Ht\ correspond to rational tensors on
\PnC[3] which may have poles of certain types along these 10 quadrics.

We shall elaborate  on this result in the subsequent lines. 
%
The map $\Ht \longto \PnC[3]$ is given by the four theta constants of the
second kind $f_0,
\dots,f_3$. These are Siegel modular forms of weight 1/2 with respect to a
common multiplier system \vf , cf. \pref{vf}.

It follows from Runge's results that the ring of all modular forms of transformation type 
$$f(\MZ)=\vf[r](M) \cdot \sqrt{\det(CZ+D)}\;\strut^r \cdot f(Z),$$
which is usually denoted by
\[ \A(\Gttf)
 =\bigoplus_{r \in \N}\left[\Gttf, r/2,\vf[r] \right], \]
is \polyfa . This theorem can be found on page \pageref{structurethm}.

The simplest case are tensors of the form \fZtail. They belong to
complex valued modular forms transforming as follows
$$f(\MZ)=\det(CZ+D)^{3\kforl} \cdot f(Z),
\hspace{2cm} M \in \Gttf.$$

\vspace{0.2cm}

Returning to the vector valued case, we start with fixing some notation.
Here and subsequently, \Mp[r] stands for the vector space of modular
forms of transformation type
$$f(\MZ)=\vf[r](M)\cdot \sqrtr{\det(CZ+D)} \cdot (CZ+D) f(Z)(CZ+D)^t,
\hspace{.3cm} M \in \Gttf.$$

It is also possible to twist this vector space by the character \vf[2]. The space
\Mm[r] consists of the modular forms satisfying
$$f(\MZ)=\vf[2](M)\cdot\vf[r](M)\cdot \sqrtr{\det(CZ+D)} \cdot (CZ+D) f(Z)(CZ+D)^t$$
for all $ M \in \Gttf$.

We shall study the graded \Attfvf-modules 
$$ \Mp(\Gttf) : = \bigoplus_{r \in \Z}\Mp[r]
\hspace{1cm}\text{and}\hspace{1cm} 
\Mm(\Gttf) : = \bigoplus_{r \in \Z}\Mm[r].
$$

The module \Mp \ contains the so called
Rankin-Cohen brackets. These are constructed by means of scalar valued
modular forms $f,g$ and derivatives, i.e.
$$\aoki{f}{g}=f\cdot Dg- g \cdot Df.$$
There is a similar construction \rcb{f}{g}{h} which defines an element
in \Mm , cf. \pref{rcbthree}.

The main results of this
thesis (\psref{structurethmeven}{structurethmodd}) are 

\vspace{.3cm}
\kasten{
\begin{align*}
\hspace{2cm}\Attfrhvf
& =\sum_{0\leq i < j \leq 3}
 \polyfab\aoki{f_i}{f_j}
\intertext{and}
\hspace{2cm}\Attfrhvfvf
& =
\sum_{
0 \leq i_1< i_2< i_3 \leq 3
}
\polyfab\rcb{f_{i_1}}{f_{i_2}}{f_{i_3}}.
\end{align*}}

For any given degree $r$, we shall exhibit explicit bases of \Mp[r] and \Mm[r].  Consequently, we obtain the
Hilbert functions 
\begin{align*}
\hspace{2cm}\dim{\Mp[r]}
= & \ 3\cdot \binom{r+1}{3}
+2\cdot \binom{r}{2}
+\binom{r-1}{1}
\intertext{and}
\hspace{2cm}\dim{\Mm[r]}
= & \ 3\cdot \binom{r-2}{3}
+\binom{r-3}{2}.
\end{align*}%

The case of modular forms with trivial multiplier system follows easily from the above lines as shown in \cref{structurethm}.

Note that the modules \Attfrhvf \  and \Attfrhvfvf \ contain the \Gttf-invariant \holten
s shown in \cref{eventype,oddtype}. Using them, we shall give a brief
overview of the structure theorems' proofs. For the sake of simplicity, we restrict ourselves to presenting the case \Attfrhvf. As shown in
\cref{algsimp}, these tensors correspond to rational tensors on~\PnC \ having poles of certain types
along 10 quadrics, that are given in \pref{thetapoly}. These tensors become \hol \ after pulling them back along
2-coverings that are ramified over the quadrics. We refer to
\pref{Omdef} for the exact wording. If the parameter~$r \in \N$ of
\cref{eventype} is even, then we can work out this condition explicitly, cf.
\pref{dimofOotcor}. This shows the equalities
$M^+_{6r}=\Mp[6r]$ for even $r$, cf.
\pref{structurethmeventen}.
The case of arbitrary~\Mp[s] can be reduced to the above
ones by multiplying with monomials in the $f_a$, cf.
\cref{structurethmeven}. This reduction uses
the very simple structure of the ring of modular forms, i.e.~$\Attfvf=\polyfa$.

This article is based on the author's PhD thesis \cite{phd} under the supervision of E.~Freitag who should be thanked for many fruitful discussions. In the thesis the proofs are presented in more detail.

\section{Symplectic groups}

The symplectic group is the subgroup 
$$\Spnr:=\set{ M \in \SL(2n,R) : \Jn\eck{M}=\Jn }$$
of the linear group of a commutative unital ring $\GL(2n,R)$ where \Jn \ denotes the involution.
The group  \SpnZ \ or  \Gn \ for short acts on the Siegel upper half-space,  i.e.  the set 
$\Hn:=\set{ Z \in \sym : \Im Z \text{ pos. \! def.} }$, by
$(M,Z) \longmapsto \MZ=(AZ+B)(CZ+D)^{-1}.$
\newline
The kernel of the natural group homomorphism \ $\SpnZ \map \SpnZq$ is the
principal congruence subgroup \Gnq. Let $A$ be a square matrix then $\diag{A}$ denotes the diagonal
vector $(a_{11},\dots , a_{nn})$.
Then Igusa's group of level $q$ is defined to be
$$\G_n\eck{q,2q}:=\set{M \in \G_n\eck{q} : \diag{AB^t}\equiv\diag{CD^t}\equiv0
\mod 2q}.$$
Suppose that $q$ is even then $\G_n\eck{q,2q}$ is normal in  the full
modular group.

\begin{lemma}\label{finitesubgroups}
The non-trivial subgroups of \Gttfpm\  of finite order are all of order
2. Their generators are conjugated in $\nicefrac{\Gt}{\pm I_4}$  to the
image of the matrix 
\[
\left(\begin{array}{cc} 
\begin{array}{cc}1 &  \\  & -1 \end{array}
 & 0 \\
0 & 
\begin{array}{cc}1 &  \\  & -1 \end{array}
 \end{array} \right) .
\]
\vspace{-0.5cm}%
\begin{proof}
It follows from rather basic algebraic
facts that the group $\G_2[4]$ acts without fixed points. Therefore, \cite{runge}'s Lemma 5.3 on page 23 completes the proof.
\end{proof}
\end{lemma}

\section{Modular~forms~of~half~integral~weight}
\label{sec:modform}
\index{modular form}%

Let \G \ be a congruence subgroup in \SpnR  ,
$v$ a multiplier system of weight $r$/2 
and \rh \ be a rational representation 
of \GLnC \ on a vector space $V$. Then,
$$v(M)\sqrtr{\det{CZ+D}}\rh(CZ+D)$$
is a factor of automorphy.
\newline
By a vector valued modular form with respect to this factor of
automorphy we mean a \holfn \ $f: \Hn \to V$ which transforms
as
$$f(\MZ)=v(M)\sqrtr{\det{CZ+D}}\rh(CZ+D)f(Z)$$
under \G.

In the case where $n=1$, the usual condition at  the cusps has to be
added. For $n\geq 2$, the Koecher principle ensures this, cf.
\cite[Hilfssatz 4.11, p.175]{frei:siegel}.

The vector space of modular forms is denoted by \vecmodfun.%
\nomenclature[\ eckmodvec]{\vecmodfun}{vector space of vector valued modular
  forms of weight $r$/2 to the representation \rh \ and the mulitplier system $v$}
If $v$ is trivial and $r=2k$ is even, then we shorten this to $\left[
  \Gamma, \left( k,\rh \right) \right]$.
\nc[\ eckmodvectriv]{$\left[
  \Gamma, \left( k,\rh \right) \right]$}{\vecmodfun\ with trivial $v$ and $r=2k$}

The pairs $(\frac{r}{2},\rh)$ and $(\frac{r-2k}{2},\det^k \cdot \rh)$
define the same factor of automorphy. Hence, for an irreducible
representation \rh \ we always may assume that \rh\ is reduced,
i.e. it is polynomial and does not vanish on $\set{\det(A)=0}$. After this
normalization we call $r/2$ the weight of a vector valued modular form.
\index{weight!modular form}

We recall  that vector valued modular forms of negative weight are identically zero.

We can construct from a multiplier system of weight 1/2 its graded algebra of
modular forms of half integral weight $A(\G,v) :=\bigoplus_{r \in
  \Z}\eck{\Gamma,\frac{r}{2},v^r}$.
  
We want to identify vector valued modular forms with tensors on \Hn ;
therefore, we fix some notation. We consider the wedge product of all $1$-forms \dZij\  in
lexicographic order by
\nc[\Altk dz]{$\dzlex$}{lexicographic ordered wedge product over all \dZij on \Hn }%
\nc[dZij]{$\dZij $}{canonical 1-form on \Hn }%
$$\dzlex=\bigwedge_{1 \leq i < j \leq n}\dZij .$$

\begin{lem}\label{modformdiffform}
There is an  isomorphism between the vector space \vecmodfunst \   of vector valued
modular forms transforming as  
$$f(\MZ)=\det(CZ+D)^{(n+1)k} \cdot (CZ+D)f(Z)(CZ+D)^t$$
and the subspace of \G -invariant tensors in $\Oot[\frac{n(n+1)}{2}](\Hn)$, i.e.%
\abb{\Phi}
{\left(\Oot[\frac{n(n+1)}{2}](\Hn)
  \right)^\G}
{\vecmodfunst}
{\om=\sum_{i<j}f_{ij}\ dZ^{ij}\ten \left(\dzlex\right)^{\ten k}}
{\left (f_{ij}\right)_{1 \leq i ,j \leq n}.}
\begin{proof}
We just have to observe the derivative of $Z \longmapsto \MZ $, cf. \cite[1.6 Hilfssatz,
p.27]{frei:siegel}.
\end{proof}
\end{lem}

We shall also consider \G -invariant tensors in $\Oiso(\Hn)$. They
can be identified with vector valued modular forms transforming as
$$f(\MZ)=\det{(CZ+D)}^{(n+1)(k+1)}(CZ+D)^{-t} f(Z) (CZ+D)^{-1}$$
under \G, cf. \cite[$4.6_1$ Folgerung, p.172]{frei:siegel}.

In the case where
$n=2$, the reduced representation~$\det^2(A) A^{-t}XA^{-1}$ is isomorphic to \rhstsym.
Therefore, we obtain the subsequent result.

\begin{lem}\label{modformtwoform}
The vector space \vecmodfunstsymt[3k+1]  is isomorphic to the subspace of \G -invariant tensors in $\Ott[3](\Ht)$.
%
The tensor
$$
\left(
f_0 \ \dZ^1\!\wedge \dZ^2 +f_1\ \dZ^0\! \wedge \dZ^2 + f_2\ \dZ^0\!\wedge
\dZ^1\right) \ten \left( \dzlex \right)^{\ten k}$$
is mapped to the vector valued modular form
$$
\begin{pmatrix}
f_2 & -f_1 \\
-f_1 & f_0
\end{pmatrix}.
$$

\end{lem}

\section{Theta series}
\label{sec:theta}

We call a vector $\mchar=(\mchar^1,\mchar^2) \in \Zn[2n]$
even if the Euclidean scalar product of $\mchar^1$ and
$\mchar^2$ 
$\inpro{\mchar^1,\mchar^2}\equiv 0 \ \mod \ 2$.

\begin{defi}[Theta series]\index{theta series!of the first kind, $\vt\eck{\mchar}$}\label{thetadef}
\index{modular form!theta series|see{theta series}}%
\index{theta series!of the second kind, $f_a$}%
On \sym \ we define  different kinds of theta series : %
\nc[theta]{$\vt\eck{\mchar}$}{theta series of the first kind}%
\begin{align*} %
\intertext{of the first kind}
\vt\eck{\mchar}(Z) := 
& \ \thesermZ  :=\sum_{g \in \Zn}{exp \left(\p i \left(%
Z\eck{g+\frac{\mchar^1}{2}}+\left(g+\frac{\mchar^1}{2}\right)^t\mchar^2%
\right)\right)} \\
& \text{\ for the characteristic $\mchar=(\mchar^1,\mchar^2) \in
\set{0,1}^{2n}\subset \Zn[2n]$; \nc[m]{\mchar}{characteristic of the
  theta series $\vt\eck{\mchar}$}}\\
\intertext{of the second kind \nc[fa]{$f_a(Z)$}{theta series of the second kind}
}
f_a(Z) :=
& \ \vt\eck{\begin{array}{c}a \\ 0 \end{array}}(2Z)
=\sum_{g \in \Zn}{exp \left(2\p i Z\eck{g+\frac{a}{2}}%
\right)}, \quad \text{for } a \in \set{0,1}^{n}.
\intertext{\parbox{\textwidth}{If desired we may observe the above parameter $a$ as an element of $(\F_2)^n$. \newline
We shall multiply the different theta series of the first kind with even characteristic}}
\theta(Z) := 
& \prod_{\mchar \in \set{0,1}^{2n} even}\vt\eck{\mchar}(Z).
\end{align*}\nc[F]{$\F_2$}{field of order 2}
\end{defi}

The theta series of the first and second kind are \holfn s on \Hn.
We state another of their properties that can be found on page 233 of \cite{igusa_theta}.

\begin{lem}\label{thetalink}
The theta series of the first and second kind are related in the following manner
$$\theserqmZ=\sum_{a \in
  (\F_2)^n}(-1)^{\inpro{a,\mchar^2}}f_{a+\mchar^1}(Z) \cdot f_{a}(Z).$$
\end{lem}

Runge deduces on page 59 of \cite{runge} that there is not a single $Z$ in \Hn[2] \  which is a common root of all four thetas, i.e. $f_a(Z)=0$. 
\newline
A fundamental result on the zero locus of Igusa's cusp form of weight 5 \chifive \ is the subsequent one, cf. Satz 2 in \cite{frei:goettingen}.

\begin{lem}\label{thetazerolocus}
The zero locus of the above \chifive\ 
on \Ht \ is $\union_{M \in \G_2}M
\inpro{\Delta_2}$, where $\Delta_2$ denotes the set of diagonal matrices.
\end{lem}

Later, we shall use the subsequent two theorems.

\begin{thm}[Structure theorem for \Attfvf]\label{structurethm}\label{vf}
The functions $f_0, \dots , f_3$ are modular forms of weight $1/2$
with respect to Igusa's group~\Gttf \ and a common multiplier system
\vf. Nota bene $\vf[4]=1$.\\
In particular, the whole ring of modular forms with respect to \Gttf \ is generated 
by the theta constants, i.e.
\nc[Ai]{\Attfvf}{graded algebra of scalar modular forms with respect to \Gttf \ and
the multiplier system \vf}
\nc[vf]{\vf}{multiplier system \vf\ induced by the theta series $f_a$,
n.b. $\vf[4]=1$}
$$\Attfvf=\bigoplus_{r \in \N}\scamodfunth[r]=\C\left[ f_0,f_1,f_2,f_3
\right].$$
\end{thm}

The modular form \chifive \ has a trivial multiplier system on \Gttf \
and hence is not contained in \Attfvf. The following variant of
the previous theorem is also true.

\begin{thm}[Structure theorem for \Attf]\label{structurethmtriv}
The graded algebra of scalar modular forms with respect to \Gttf \ and
the trivial multiplier system can be decomposed
\nc[Aint]{\Attf}{graded algebra of scalar modular forms with respect to \Gttf \ and
the trivial multiplier system}
\begin{align*}
\Attf
& =\bigoplus_{r \in \N}\left[\Gttf, r \right]\\
& =\left(\bigoplus_{d\geq 0} \C_{4d}\left[ f_0,f_1,f_2,f_3\right]\right)
\ \bigoplus\ \left(\bigoplus_{d\geq 0} \C_{4d}\left[ f_0,f_1,f_2,f_3\right]\right) \cdot \chifive
\end{align*}
\end{thm}

This version can be found in \cite[Remark 3.15, p.74]{runge}. And
\cref{structurethm} is a rather simple consequence using the methods
of \cref{sec:modular_existence_even,sec:modular_existence_odd,sec:modular_existence_int,sec:sec:modular_existence_scalar}, cf.  \pref{structurethmAttfvf}.

\section{Covering \hol \ tensors}
\label{sec:covhol}


We call a \hol, surjective and finite covering map $f:M \map N$
between two $n$ dimensional complex manifolds \entcov \ map
if \rami{f}{M} is a smooth hypersurface.  Prototypes are the $k$-th $n$-dimensional standard elements 
\abb{p_n^k }{\E^n}{\E^n}{(z^1, \dots , z^n)}{(z^1, \dots ,
  z^{n-1},\left( z^n \right)^k),}%
 where $n,k>0$. Indeed, any simple covering is locally isomorphic to a standard element.
This follows from the topological classification of unramified
covering maps of the space~$\En[n-1]\times\E^*$. 
\newline
Let $R$ be a closed submanifold of codimension 1 of $N$. Then, for each
point $p \in R$ there is a chart $V \to \En$ that sends $R$ to
$z^n=0$. Hence, there is a simple covering $U \to V$ that ramifies 
over $R \cap V$.

\begin{defi}[Covering \hol \ tensors]\label{Omdef}\index{covering!covering \hol \ tensor}
\index{tensor!holomorphic!covering \hol \ tensor}
Let $D$ be an effective divisor on a $n$-dimensional complex manifold $M$, we define
$\Ot(M,D)$\nc[OZaz]{$\Ot(M,D)$}{vector space of covering \hol \ tensors}
 as the space of tensors $\om \in \Ot(M\ohne \supp D)$ with
the supplementary property :%
\begin{addmargin}[1cm]{1cm}
Let $Y$ be an irreducible component of $D$ and $q$ a point in $Y$  that is
 a regular point in $\supp D$. Then, there exist an open \nhd \ of
$q\in V$ and \anentcov \ $p:U \to V$ with the subsequent properties :
\begin{enumerate}\setlength{\itemsep}{-8pt} 
\item
$p$  is ramified over $Y \cap V$;
\item
$p$ is isomorphic to the standard element $\pnk$, where $k=D(Y)+1$;
\item
$\om$'s pullback $p^*\om$ is \hol ally extendable to the whole of $U$.
\end{enumerate}
\end{addmargin} 
\end{defi}
Of course, the third condition in the previous defintion is independent of the chosen $p$.
It is not hard to show that these tensors are meromorphic on $M$.
Later, we shall use tensors of the following type (the notation should be self explicatory) 
$$\left( \Omp \ten (\Omn )^{\ten k} \right) (M,D), \quad \text{where }n=\dim{M}.$$
\nc[OZazpk]{$(\Omn )^{\ten k} (M,D)$}{vector
  space of covering \hol \ tensors of type $(\Omn )^{\ten k}$}


\section{The quotient space
  \texorpdfstring{$\lecos{\ \Hn[2]}{\Gttf}$}{(Gamma_2 [2,4])\textbackslash IH_2}}

The considerations in \cref{sec:covhol} were aimed at the group \Gttf \ acting on \Ht \ and the induced quotient manifold $\lecos{\Hn[2]}{\Gttf}$. \Pref{finitesubgroups} implies
that all elements of finite order in \Gttf\ are conjugated to
 the matrix
\[
\left(\begin{array}{cc} 
\begin{array}{cc}1 &  \\  & -1 \end{array}
 & 0 \\
0 & 
\begin{array}{cc}1 &  \\  & -1 \end{array}
 \end{array} \right).
\]
The fixed point set of \Gttf\  is the disjoint union of the \Gn[2] conjugates of the diagonal~\Dn[2]. Its image~\ramdiv\  in $\lecos{\Hn[2]}{\Gttf}$ is the ramification divisor (its multiplicities are just 1). We deduce from \pref{thetazerolocus}  the succeeding lemma.

\begin{lem}\label{Szerotheta}
The ramificatin divisor \ramdiv \ is 
the zero locus of $\chifive =\prod_{\mchar \in \set{0,1}^4 even}
\vt\eck{\mchar}$. 
Moreover, the image of the zero locus of an individual
$\vt\eck{\mchar}$ is one of the above mentioned 10 components
of the type $\MZ[{\Dn[2]}]$.
\end{lem}

The four theta series of the second kind $f_0,\dots, f_3$ have no common root. Therefore they induce a \hol \ map
$$
F: \Hn[2] \longmap \PnC[3], \quad
   Z \longmapsto \eck{f_0(Z), \dots ,f_3(Z)}.
$$

From Runge's results \cite{runge} one can deduce the following
striking theorem.

\begin{thm}
The map $\phi : \lecos{\Ht}{\Gttf} \to \PnC[3]$ is an open,
holomorphic embedding. The image's complement is
an analytic subset of codimension 2.
\end{thm}

Instead of observing divisors, \fn s, and tensors on \lecos{\Ht}{\Gttf}
we can study their counterparts on \PnC[3].
In particular, we can consider \ramdiv \ as a  divisor 
 on \PnC[3].

\begin{lem}\label{thetapoly}
The ramification divisor \ramdiv\ considered on \PnC[3] is the sum of ~10 quadrics given by the~10 polynomials known from \pref{thetalink}.
Prototypes are  $Q_1 (z^0,\dots, z^3 )= (z^0)^2 - (z^1)^2 + (z^2)^2 - (z^3)^2$ and
$Q_8 (z^0,\dots, z^3 )= 2(z^0z^2 + z^1z^3)$.
\end{lem}

The map $\Ht \longmap \lecos{\Ht}{\Gttf}$
is locally a simple covering ramified over \ramdiv . Therefore we obtain the following result.

\begin{thm}\label{algsimp}
There is a natural isomorphism
$$\Om^{\ten q} (\Ht)^\Gttf \iso \Om^{\ten q}(\PnC[3],\ramdiv).$$
\end{thm}

The tensors on the left hand side can be considered as certain vector
valued Siegel modular forms and the ones on the right hand side can be
easily described in algebraic terms.

\label{chap:main}


\section{Construction of \meroten s with prescribed poles}

For the construction of these tensors, we consider a homogeneous polynomial $Q$ of degree~$d$ in $n+1$ variables
$X^0, \dots , X^{n}$. We always assume that $Q$ is square-free and
that~$X^0$ does not divide~$Q$.
We want to describe \meroten s on \PnC \ that
are \hol \ outside~$Q$'s zero locus \zero{Q}. 
Recall  that the
projective coordinates of \PnC \ are denoted by $z^0, \dots ,z^n$ and
the coordinates on the affine space \Anull\ are 
$$ \left(\zt^1, \dots ,\zt^n\right):=\left(\frac{z^1}{z^0}, \dots , \frac{z^n}{z^0}\right).$$ 
For the sake of simplicity we take \om \ to be of the type \ \Opt.
We introduce a handy notation for the canonical basis elements of \Omp.
Given a subset 
$I \subset \set{1,\dots, n}$, say
$I=\set{i_1,\dots, i_p} $  
and  
$i_1 < \dots < i_p$,
then \dzti[I] \ is short for 
$$\dzti[I]=\dzti[i_1]\we \dots \we \dzti[i_p].$$

It is not hard to show the following lemma.
 
\begin{lem}\label{OitPnCQvectors}
Let \om \ be a \meroten \ on \PnC\  which is \hol \ outside \zero{Q},
then \om \ can be written in the form 
$$\om=\sum_{
\substack{I \subset \set{1,\dots, n}\\
 \abs{I}=p}
}\om_{I} \ \dzti[I] \ten \left( \dzti[1]\wedge
\dots \wedge\dzti[n]\right)^{\ten k},$$
where
$$ \om_{I}(\zt^1, \dots, \zt^n)=\frac{A_{I}(1,\zt^1, \dots, \zt^n)}{Q^N(1,\zt^1, \dots, \zt^n)}$$
with the following properties \vspace{-0.2cm}
\begin{enumerate} \setlength{\itemsep}{-6pt} 
\item $N$ is a natural number;
\item $A_{I}$
is a homogeneous polynomial of degree $N \cdot \deg{Q}$;
\item $(X^0)^{k(n+1)+p}$ divides $A_I$;
\item and it holds for all $J \subset \set{2,\dots,n}$ with
  $\abs{J}=p-1$
 $$\left. (X^0)^{k(n+1)+p+1}\ \right| \ \sum_{
\substack{1 \leq  j \leq n \\
 j \notin J }} (-1)^{pos(j,J \cup \set{j})}\cdot X^j \cdot
A_{J\cup\set{j}} ,$$
where $pos(i,I)$ returns $i$'s position in the ordered set $I$. 
 \end{enumerate} 
For $p=n$ conditions 3 and 4 are merged to  
\begin{enumerate}
\item[3'.] $(X^0)^{k(n+1)+n+1}$ 
divides $A_{\set{1, \dots, n}}$.
\end{enumerate}
\end{lem}

We want to describe the space $ \Opt (\PnC,\ramdiv)$  for the divisor
$\ramdiv=\sum \am \zero{\Qm}$ where $\am\geq 0$.

\begin{thm}\label{dimofOot}
We introduce the numbers
\footnote{in the following formula \floor{x} denotes the floor
  function $\max\, \{m\in\mathbb{Z}\mid m\leq x\}$.} 
\nc[\ ceil]{$\floor{x}$}{floor function of the real number $x$}
$$d:=\max_{1\leq \mchar \leq \mcharanz}\floor{\frac{\am}{\am +1}k}
$$ and 
$$D:=\max_{1\leq \mchar \leq \mcharanz}\floor{\frac{\am}{\am
    +1}(k+1)}.$$

Let \om \ be a tensor as in \cref{OitPnCQvectors} where we choose $N$
minimal.
We can state \nec \ and sufficient conditions in terms of $d$ and $D$
for \om \ to lie in $ \Opt (\PnC,\ramdiv)$~:
$$\om \in \Opt (\PnC,\ramdiv) \quad \longhence \quad N \leq D,$$
and
$$N \leq d \quad \longhence \quad  \om \in \Opt (\PnC,\ramdiv).$$
%
%
\begin{cor}\label{dimofOotcor}
If it holds for all $\mchar$
$$\am +1 \ \vert\ k, $$
then $d$=$D$. Hence, the space $\Opt (\PnC,\ramdiv)$ is completely determined.
\end{cor}
\begin{proof}[Theorem]
We shall just present the proof in the case where each multiplicity \am\ and the exterior power
 $p$  equal 1.
\newline
Let $x$ be a smooth point of $Q$'s zero locus and \Qm \  its prime
factor that vanishes at $x$. We have to consider an open \nhd \ $x \in
V \subset \PnC$ and \anentcov \ $U \to V$ which ramifies over $\Qm=0$.
 \Obda, we may assume that $V$ is a subset of
\Anull\ and that the $n$-th partial derivative of \Qm \ is
non-zero in $V$.
Now, we can use the considerations from the beginning of \cref{sec:covhol} for the
submanifold $V \cap \zero{\Qm}$. The covering $p$ is of the form 
$$(\zt^1,\dots,\zt^n)\longmapsto (\zt^1,\dots , \zt^{n-1},(\zt^n)^{2}+\varphi(\zt^1, \dots , \zt^{n-1})),$$
where $\varphi$ is implicitly defined by
$$\Qm(1,\zt^1,\dots, \zt^n)=0 \quad \iff \quad \zt^n=\varphi(\zt^1,\dots, \zt^{n-1}).$$
We have to observe the pullback $p^*\om$ of the tensor
$$\om=\sum_{i+1}^n\om_{i} \ \dzti \ten \left( \dzti[1]\wedge
\dots \wedge\dzti[n]\right)^{\ten k}.$$
The coefficient functions are
\begin{align*}
\label{pullbackcoordinates}
(p^*\om)_n & = \left(2 \cdot \zt^n\right)^{k+1}\  \om_n
\circ p 
\text{\hspace{8.4em}  
and}\\
 (p^*\om)_j & =\left(2 \cdot \zt^n\right)^k 
\left( \om_j \circ p+\om_n \circ p \cdot \delztj{\varphi}\right)
\text{\hspace{2em}  
for $j\neq n$.} 
\end{align*}
We have to check whether the pullback $p^*\om$ is \hol \ on $U$, i.e.
$$\ord{(p^*\om)_j,\zero{\zt^n}}\geq 0, \quad 1\leq j \leq n.$$
In the case where $j=n$, this is equivalent to
$$\ord{\om_n\circ p,\zero{\zt^n}}\geq -(k+1).$$
The ultrametric inequality yields that it is necessary for $ p^*\om$ to be \hol\  that all coefficients  satisfy the above inequality. In the same way we derive sufficient conditions for $p^*\om$ to be \hol . They are
$$\ord{\om_j\circ p,\zero{\zt^n}}\geq - k \quad \forall\  j \in
\set{1,\dots n} .$$
Now, we want to consider the coefficient functions $\om_j$ instead
of the functions $\om_j~\circ~p$. The zero set \zero{\zt^n} in $U$
corresponds to the zero divisor $(\Qm)$ of the polynomial~$\Qm(1,\zt^1,\dots, \zt^n)$ on $V$.
We have
$$\ord{\om_j\circ p,\zero{\zt^n}}=2 \cdot \ord{\om_j,(\Qm)}.$$
So, we obtain the necessary condition 
$$\ord{\om_j,(\Qm)} \geq - \floor{\frac{k+1}{2}}, \quad
1\leq j \leq n,$$
and the sufficient condition
$$\ord{\om_j,(\Qm)} \geq - \floor{\frac{k}{2}}, \quad
1\leq j \leq n,$$
for \om \ to be in the space
 $\Oot(\PnC,\ramdiv)$.
This completes the proof of the theorem.
\end{proof}
\end{thm}

\section[A structure theorem for vector valued modular forms w.r.t. the multiplier system \vf]{A structure theorem for vector valued modular forms with respect to the multiplier system \texorpdfstring{$\mathbf \vf$}{v_f}}
\label{sec:modular_existence_even}

We consider the polynomial $Q=\prod_\mchar \Qm $ of degree 20 that gives the divisor of
\chifive\  in \PnC[3], cf. \pref{thetapoly}. 
We recall that 
$$\left[ \Gttf,\left( \frac{12k}{2},\rhstsym \right),\vf[12k] \right]
\iso \left(\Om \tenO \left(\Omp[3]\right)^{\ten 2k}\right)
(\PnC[3],(Q))$$
according to \prefs{modformdiffform}{algsimp}. Hence, we can reformulate the results of  the previous section in
terms of modular forms.

\begin{lem}\label{detequalschi}
If we denote by $f_a$ the theta constants of the second kind, then the determinant 
$$
\begin{vmatrix}
\frac{\del \left(\frac{f_{1}}{f_0}\right)}
{\del Z^{0} } 
& \frac{\del\left(\frac{f_{1}}{f_0}\right)}
{\del Z^{1} } 
& \frac{\del\left(\frac{f_{1}}{f_0}\right)}
{\del Z^{2} } 
\\[0.3 em]
 \frac{\del\left(\frac{f_{2}}{f_0}\right)}
{\del Z^{0}}
 &  \frac{\del\left(\frac{f_{2}}{f_0}\right)}
{\del Z^{1}}
 &  \frac{\del\left(\frac{f_{2}}{f_0}\right)}
{\del Z^{2}}
\\[0.3 em]
\frac{\del \left(\frac{f_{3}}{f_0}\right)}
{\del Z^{0} } 
& \frac{\del\left(\frac{f_{3}}{f_0}\right)}
{\del Z^{1} } 
& \frac{\del\left(\frac{f_{3}}{f_0}\right)}
{\del Z^{2} } 
\end{vmatrix}
$$
equals $$c_5 \cdot \frac{\chifive}{(f_0)^4}$$ 
with a constant $c_5$ in \C .
\begin{proof}
A proof can be found on pages 15 and 16 of \cite{FS}.
\end{proof}
\end{lem}

We shall present an easy and  well-known lemma.

\begin{lem}\label{rcbtwo}
For $f$ and $g$ in \scamodfunth \ the
Rankin-Cohen bracket
\nc[\ rcb]{$\aoki { f}{ g}$}{Rankin-Cohen bracket of the modular forms
$f$ and $g$}
\newline
\begin{align*}
\aoki { f}{ g}
& :=
f\cdot Dg- g \cdot Df
=
\begin{pmatrix}
\left( f\frac{\del g}{\del Z^0} 
-g\frac{\del f}{\del Z^0} \right)
 &
\left( f\frac{\del g}{\del Z^1} 
-g\frac{\del f}{\del Z^1} \right)
\\
\left( f\frac{\del g}{\del Z^1} 
-g\frac{\del f}{\del Z^1} \right)
&
\left( f\frac{\del g}{\del Z^2} 
-g\frac{\del f}{\del Z^2} \right)
\\
\end{pmatrix}
=f^2 \cdot D\left( \frac{g}{f}\right)
\end{align*}
%
%
%
lies in $\eck{\Gttf,(1,\rhstsym),\vf^2}$.
\end{lem}

The next lemma will rely heavily on the just defined Rankin-Cohen brackets.

\begin{lem}\label{dimofvecmodzero}
Every modular form $f \in \vecmodfunstpm[6\kforls]$ is of the form
$$f=\sum_{1\leq i \leq 3} P_i(f_0,\dots f_3) \aoki{f_0}{f_i} \inv{f_0},$$
where all $P_i$ are homogeneous polynomials of degree $12\kforls-1$.
Conversely, such a sum lies in \vecmodfunstpm[6\kforls] iff it is \hol
\ which means 
$$f_0 \ \vert \sum_{1\leq i \leq 3} f_i  \cdot P_i .$$
\begin{proof}
A modular form 
$f \in \vecmodfunstpm[6\kforls]$
 can be considered as tensor on \PnC[3] of the form
$$\sum_{i=1}^3 R_i\  \diff \left(\frac{f_i}{f_0} \right)\ten 
\left(\diff\left(\frac{f_1}{f_0}\right) \wedge \cdots \we \diff \left(\frac{f_3}{f_0} \right) \right)^{\ten 2\kforls }.$$
In \pref{dimofOot}, we have seen that $R_i$ is of the form
$$R_i=\frac{P_i \cdot f_0^{8 \kforls+1}}{\chifive^{2\kforls}}$$
with $P_i$ a polynomial of degree $12 \kforls-1$ and
$$f_0 \ \vert \sum_{1\leq i \leq 3} f_i  \cdot P_i .$$
We want to observe this tensor on the upper half plane. The functional
determinant of~$\diff\left(\frac{f_1}{f_0}\right) \wedge \cdots \we \diff
\left(\frac{f_3}{f_0} \right)$ 
is $\chifive / f_0^4$, cf. \pref{detequalschi}.
Hence, the modular form is of the desired type.
\end{proof}
\end{lem}

We could have formulated the above theorem replacing $f_0$ by any other $f_a$.

We introduce the \polyfa -module
\nc[Mmod]{${\cal M}^+(\Gttf)$}{\polyfa -module of vector valued modular
  forms with respect to \Gttf, the multiplier system \vf\ and \rhsym}
$${\cal M}^+(\Gttf):=\bigoplus_{r \in \Z}\vecmodfunstvf. $$
One of its  \polyfa-submodules is
\nc[Mmod]{$M^+$}{\polyfa -module of Rankin-Cohen brackets \aoki{f_i}{f_j}}
$$M^+:=\sum_{0\leq i < j \leq 3}
 \polyfa\aoki{f_i}{f_j}.$$

\begin{thm}\label{aokirelation}
The Rankin-Cohen brackets are related in the following manner :
\begin{align*}
R_1 : f_1\aoki{f_0}{f_2}
& =f_2 \aoki{f_0}{f_1} + f_0\aoki{f_1}{f_2}\ ;%
\\
R_2 : f_1\aoki{f_0}{f_3}
& = f_3 \aoki{f_0}{f_1} +f_0\aoki{f_1}{f_3} \ ;
\\
R_3 : f_2\aoki{f_0}{f_3}
& =f_3 \aoki{f_0}{f_2} + f_0\aoki{f_2}{f_3}\ ; %
\\
R_4 : f_2\aoki{f_1}{f_3}
& =f_3 \aoki{f_1}{f_2} + f_1\aoki{f_2}{f_3}\ .
\end{align*}
These are defining relations of the module $M^+$. Therefore, the Hilbert function
is $$
\dim{M^+_r} 
= 3\cdot \binom{r+1}{3}
+2\cdot \binom{r}{2}
+\binom{r-1}{1}
.$$
\begin{proof}
Let $R$ be an arbitrary relation 
$$\quad \sum_{0 \leq i<j\leq 3} P_{ij}\aoki{f_i}{f_j}=0.$$
It is equivalent to the subsequent three relations
$$f_0P_{0j}
+\sum_{i=1}^{j-1}f_iP_{ij}
-\sum_{i=j+1}^3 f_iP_{ji}=0, \quad
\forall \ j \in \set{1,\dots,3}$$
due to the given relations $R_1, \dots, R_4$ and the linear independence of the brackets  \aoki{f_0}{f_j}.

Applying $R_1, \dots, R_4$, the relation  $R$  can be transformed to a form where $P_{ij}$ is a polynomial in the
variables $f_0, \dots, f_j$ .
In this normal form we see that each $P_{ij}$ is zero. Indeed, for
$j=1$ we obtain
$$f_0P_{01}-\sum_{i>1}f_iP_{1i}=0.$$
Setting the variables $f_2, f_3$ zero yields
$$P_{01}=P_{01}(f_1)=0.$$
Specialising now $f_3$ gives $P_{12}=P_{12}(f_1,f_2)=0$
and hence $P_{13}=0$. 
\\
The just proven equality $P_{12}=0$ simplifies the above relation for
$j=2$ to
$$f_0P_{02}-f_3P_{23}=0. $$
A similiar line of argument shows
$$P_{02}=P_{23}=0. $$
Setting $j=3$ gives the remaining coefficients.
\end{proof}
\end{thm}

The module $M^+$ can be considered as a submodule of the free module
$$\sum_a \polyfa  \cdot Df_a.$$
In this setting, an element of $M^+$, say $\sum_a P_a \cdot Df_a$, can be characterized by
a single \polyfa - linear equation :
$$\sum_a f_a P_a =0 .$$ 

A short \magma \ code or the above appealing observation, that is due to U.
Weselmann, both imply the subsequent theorem.

\begin{thm} \label{intersectioneven}
We have
$$\bigcap_{i=0}^3 \frac{f_0 \cdots f_3 }{f_i} M^+  = f_0 \cdots f_3 \cdot M^+.$$
\end{thm}

An immediate consequence is the succeeding lemma.

\begin{lem}\label{structurethmeventen}
Every modular form $f \in \vecmodfunstpm[6\kforls]$ is of the form
$$f=\sum_{0\leq i < j \leq 3} P_{ij}(f_0,\dots f_3) \aoki{f_i}{f_j}, $$
where any $P_{ij}$ is a homogeneous polynomials of degree $12\kforls-1$.
\end{lem}

We want to generalize this now to arbitrary weights.

\begin{thm}[Structure theorem]\label{structurethmeven}
We have
$${\cal M}^+(\Gttf)=\bigoplus_{r \in \Z}\vecmodfunstvf =\sum_{0\leq i < j \leq 3}
 \polyfa\aoki{f_i}{f_j}.$$
\begin{proof}
Let $f$ be a modular form of weight $r/2$.
If $r$ is a multiply of 12, then we apply \cref{structurethmeventen}. Therefore, we may assume that $f_i \cdot f$ lies in the right
hand side. Now, we can apply \cref{intersectioneven}.
\end{proof}
\end{thm}

\section{A structure theorem for vector valued modular forms twisted
  with  the character 
  \texorpdfstring{$\mathbf \vf^2$}{v_f^2}}
\label{sec:modular_existence_odd}

We consider the character $\vf^2$ on \Gttf \ and the twisted version
of $\cal M^+$ :
$${\cal M^-}(\Gttf):=\bigoplus_{r \in \Z}\vecmodfunstvfvf.$$
It is possible to show that \Mm \ coincides with its \polyfa -submodule
$$\Mmrcb:=\sum_{0\leq i < j < k \leq 3}
 \polyfa\rcb{f_i}{f_j}{f_k}.$$%
For this purpose, we shall only present the proclaims that differ
substantially from their counterparts in
\cref{sec:modular_existence_even}.

\begin{defi}\label{rcbthree}
For $f$, $g$ and $h$ in \OU[\Ht] \ 
and $0 \leq j_1 < j_2 \leq 2$
we define 
\[
\rcb[(j_1, j_2)]{ f}{g }{h}
:=
\begin{vmatrix}
\frac{\del f}{\del Z^{j_1}} 
&
\frac{\del f}{\del Z^{j_2}} 
& 
f
\\[0.3em]
\frac{\del g}{\del Z^{j_1}} 
&
\frac{\del g}{\del Z^{j_2}} 
&
g \\[0.3em]
\frac{\del h}{\del Z^{j_1}}
&
\frac{\del h}{\del Z^{j_2}}
&
 h \\
\end{vmatrix}
.\]
Then the Rankin-Cohen 3-bracket of $f$, $g$ and $h$ is
\[
\rcb{ f}{g }{h}
:=
\begin{pmatrix}
\textcolor{white}{-}\rcb[(1, 2)]{ f}{g }{h} & -\rcb[(0, 2)]{ f}{g }{h} \\[0.3em]
-\rcb[(0, 2)]{ f}{g }{h} & \textcolor{white}{-}\rcb[(0, 1)]{ f}{g }{h} 
\end{pmatrix}.
\]%
\end{defi}%
A direct computation shows that it holds 
$$\rcb{f}{g}{h} \in \vecmodfunstvfvf[5]$$ 
for all modular forms $f,g,h$ in $[\Gttf, \frac{1}{2},\vf]$. 
Similarly as in \pref{dimofvecmodzero}, we can relate the 2-form $d\left(\frac{f_1}{f_0}\right)\we d\left(\frac{f_2}{f_0}\right)$ to the 3-bracket  \rcb{f_0}{f_1}{f_2} using \prefs{modformtwoform}{algsimp}.
In contrast to \pref{aokirelation}, the Rankin-Cohen 3-brackets just satisfy a
single relation.
\begin{lem}\label{rcbrelation}
The Rankin-Cohen 3-brackets are related in the following manner
$$
R_5 : -f_0 \rcb{f_1}{f_2}{f_3}
+f_1 \rcb{f_0}{f_2}{f_3}
-f_2 \rcb{f_0}{f_1}{f_3}
+f_3 \rcb{f_0}{f_1}{f_2}
=0.
$$
This is a defining relation of the module $\Mmrcb$. Hence, the Hilbert function
is $$
\dim{\Mmrcb[r]} 
= 3\cdot \binom{r-2}{3}
+\binom{r-3}{2}
.$$
\begin{proof}
We just show the relation. Consider the matrix
\begin{align*}
A & =\begin{pmatrix}
f_0
& f_1
& f_2
& f_3
\\[0.3 em]
\frac{\del f_0}
{\del Z^{0} } 
& \frac{\del f_1}
{\del Z^{0} } 
& \frac{\del f_2}
{\del Z^{0} } 
& \frac{\del f_3}
{\del Z^{0} } 
\\[0.3 em]
\frac{\del f_0}
{\del Z^{1} } 
& \frac{\del f_1}
{\del Z^{1} } 
& \frac{\del f_2}
{\del Z^{1} } 
& \frac{\del f_3}
{\del Z^{1} } 
\\[0.3 em]
\frac{\del f_0}
{\del Z^{2} } 
& \frac{\del f_1}
{\del Z^{2} } 
& \frac{\del f_2}
{\del Z^{2} } 
& \frac{\del f_3}
{\del Z^{2} } 
\end{pmatrix}
\intertext{and its adjoint matrix}
\Adj(A) & =\begin{pmatrix}
\frac{\del f_{\set{1,2,3}}}{\del Z^{\set{0,1,2}}}
& -\rcb[(1,2)]{f_1}{f_2}{f_3}
& +\rcb[(0,2)]{f_1}{f_2}{f_3}
& -\rcb[(0,1)]{f_1}{f_2}{f_3}
\\[0.3 em]
\vdots
& \vdots & \vdots & \vdots
\\
-\frac{\del f_{\set{0,1,2}}}{\del Z^{\set{0,1,2}}}
& +\rcb[(1,2)]{f_0}{f_1}{f_2}
& -\rcb[(0,2)]{f_0}{f_1}{f_2}
& +\rcb[(0,1)]{f_0}{f_1}{f_2}
\end{pmatrix}.
\end{align*}
The product $A\cdot\Adj(A)$ is just $c_5 \cdot \chifive \cdot I_4$;
this yields the desired relation.
\end{proof}
\end{lem}

The module \Mmrcb \ can be considered as a submodule of the free module
$$\sum_{0 \leq a < b \leq 3} \polyfa  \cdot Df_a \we Df_b.$$
In this setting, an element of \Mmrcb , say $\sum_{0 \leq a < b \leq 3} P_{ab}  \cdot Df_a \we Df_b$, can be characterized by
four \polyfa - linear equations :
\begin{align*}
f_1 P_{01} + f_2 P_{02} +f_3 P_{03}=0 ;\\
-f_0 P_{01} + f_2 P_{12} +f_3 P_{13}=0 ; \\
-f_0 P_{02}  - f_2 P_{12} +f_3 P_{23}=0 ; \\
-f_0 P_{03} - f_1 P_{13} -f_2 P_{23}=0 .
\end{align*}

This leads eventually to the aforementioned result.

\begin{thm}[Structure theorem]\label{structurethmodd}
We have
$${\cal M}^-(\Gttf)=\bigoplus_{r \in \Z}\vecmodfunstvfvf =\sum_{0\leq
  i < j < k \leq 3}
 \polyfab\rcb{f_i}{f_j}{f_k}.$$
\end{thm}

\section[A structure theorem for vector valued modular forms w.r.t. the
  trivial mult. system]{A structure theorem for vector valued modular forms with
  respect to the trivial multiplier system}
\label{sec:structurethmint}

\label{sec:modular_existence_int}

We can extract from the modules ${\cal M}^+(\Gttf)$ and ${\cal M}^-(\Gttf)$ the modular
forms with trivial multiplier system.

\begin{thm}\label{dimofvecmod}
We have
\begin{align*}
& \bigoplus_{k \in \Z}\vecmodfunstpm[k] = \\
& \sum_{0\leq i < j \leq 3}
 \polyfab[2+4\Z]\aoki{f_i}{f_j} 
\oplus \sum_{
0 \leq i_1< i_2< i_3 \leq 3
}
\polyfab[1+4\Z]\rcb{f_{i_1}}{f_{i_2}}{f_{i_3}}
\end{align*}
and
$$
\dim{\vecmodfunstpm[k]}
= 
\begin{cases}
3\cdot \binom{2k+1}{3}
+2\cdot \binom{2k}{2}
+\binom{2k-1}{1}, &
\text{\ \ \ if $k$ is even,} \\
3\cdot \binom{2k-2}{3}
+\binom{2k-3}{2},
& \text{\ \ \ if $k$ is odd.} \\
\end{cases}
$$
\end{thm}

\section{A structure theorem for scalar valued modular forms}
\label{sec:sec:modular_existence_scalar}
So far, we have treated tensors of the type
$$\Opt[3],$$
where $p=1,2$. It is worthwhile to mention that our method is also
successful for $p=3$. In this case we get the structure theorems for
the ring of scalar valued modular forms. Recall that our method relied
on the injectivity and ramification behaviour of the map
$$\eck{f_0,\dots,f_3}: \lecos{\Ht}{\Gttf} \longto \PnC[3].$$

More precisely, we obtain the following theorem which is essentially
due to Runge.

\begin{thm}\label{structurethmAttfvf}
We have
\begin{align*}
\bigoplus_{r \in \Z}\scamodfunth[r] & =\polyfa . \\
\intertext{Twisting with \vf[2]  yields}
\bigoplus_{r \in \Z}\scamodfunthvf[r] & =\chifive \cdot \polyfa . \\
\intertext{As a consequence,}
\bigoplus_{r \in \Z}\eck{\Gttf,r} & =\polyfa[4\Z] \oplus
\polyfa[4\Z] \cdot \chifive. 
\nc[\ eckmodsca]{\eck{\G,k}}{vector space of scalar valued modular forms
  of weight $k$}
\end{align*}
\end{thm}
We omit the details of the proof, but stress that we start with tensorial
weights, i.e. $k \in 3\Z$, as in the case of vector valued modular
forms. Afterwards, we extend the results to all weights again by
intersecting appropriate modules.
It is quite interesting 
that  Igusa's modular form \chifive \ comes up automatically in
our approach. This happens while studying the holomorphicity of tensors
by means of the ramification behaviour of~\eck{f_0,\dots,f_3}.

\bibliographystyle{halpha}
\bibliography{bookscited}

\begin{tabular}{l@{\hspace{1.5cm}}r}
Thomas Wieber, & \href{mailto:Thomas.Wieber@mathi.uni-heidelberg.de}{\nolinkurl{Thomas.Wieber@mathi.uni-heidelberg.de}}\\
Mathematisches Institut, & \\
Im Neuenheimer Feld 288, & \\
69120 Heidelberg, Germany & \href{http://www.mathi.uni-heidelberg.de/~twieber}{\nolinkurl{http://www.mathi.uni-heidelberg.de/~twieber}} 
\end{tabular}

\end{document}